\documentclass[12pt]{article}
\usepackage{amsmath,latexsym,epsf}
\begin{document}

\newcommand{\nc}{\newcommand}
\def\PP#1#2#3{{\mathrm{Pres}}^{#1}_{#2}{#3}\setcounter{equation}{0}}
\def\ns{$n$-star}\setcounter{equation}{0}
\def\nt{$n$-tilting}\setcounter{equation}{0}
\def\Ht#1#2#3{{{\mathrm{Hom}}_{#1}({#2},{#3})}\setcounter{equation}{0}}
\def\qp#1{{${(#1)}$-quasi-projective}\setcounter{equation}{0}}
\def\mr#1{{{\mathrm{#1}}}\setcounter{equation}{0}}
\def\mh#1{{{\mathcal{#1}}}\setcounter{equation}{0}}
\newtheorem{Th}{Theorem}[section]
\newtheorem{Def}[Th]{Definition}
\newtheorem{Lem}[Th]{Lemma}
\newtheorem{Pro}[Th]{Proposition}
\newtheorem{Cor}[Th]{Corollary}
\newtheorem{Rem}[Th]{Remark}
\newtheorem{Exm}[Th]{Example}
\def\Pf#1{{\noindent\bf Proof}.\setcounter{equation}{0}}
\def\>#1{{ $\Rightarrow$ }\setcounter{equation}{0}}
\def\<>#1{{ $\Leftrightarrow$ }\setcounter{equation}{0}}
\def\bskip#1{{ \vskip 20pt }\setcounter{equation}{0}}
\def\sskip#1{{ \vskip 5pt }\setcounter{equation}{0}}
\def\midskip#1{{ \vskip 10pt }\setcounter{equation}{0}}
\def\bg#1{\begin{#1}\setcounter{equation}{0}}
\def\ed#1{\end{#1}\setcounter{equation}{0}}
\def\r#1{{\rm{#1}}}


\title{\bf Finitistic and Representation Dimensions}
\smallskip
\author{\small {Jiaqun WEI} {\thanks {Supported by the National Science Foundation of China
(No.10601024)}}\\
\small Department of Mathematics,
Nanjing Normal University \\
\small Nanjing 210097, P.R.China\\ \small Email:
weijiaqun@njnu.edu.cn}
\date{}
\maketitle
\baselineskip 18pt
%
%
\begin{abstract}
\vskip 10pt%

Let $A$ be an artin algebra with representation dimension not more
than 3. Assume that $_AV$ is an Auslander generator and
$M\in\mr{add}_AV$ (for example,
$M\in\mr{add}_A(A\oplus\mathbf{D}A)$), we indicate that both
$\mr{findim}(\mr{End}_AM)$ and $\mr{findim}(\mr{End}_AM)^{op}$ are
finite and consequently, the Gorenstein Symmetry Conjecture, the
Wakamatsu-tilting conjecture and the generalized Nakayama
conjecture hold for $\mr{End}_AM$. In particular, we see that the
finitistic dimension conjecture and all the above conjectures hold
for artin algebras which can be realized as endomorphism algebras
of modules over representation-finite algebras. It is also shown
that if every quasi-hereditary algebra has a left idealized
extension which is a monomial algebra or an algebra whose
representation dimension is not more than 3, then the finitistic
dimension conjecture holds.
\midskip\
{\it MSC:}\ \ 16E05; 16E10; 16G20
%
%

{\it Key words:}\ \ finitistic dimension; representation
dimension; representation-finite algebra
\end{abstract}
%
\vskip 30pt
%

\section{Introduction}
 \hskip 15pt

Throughout the paper, we work on artin algebras and finitely
generated left modules.

Let $A$ be an artin algebra. Recall that the little finitistic
dimension of $A$, denoted by $\mr{findim}A$, is defined to be the
supremum of the projective dimensions of all finitely generated
modules of finite projective dimension. Similarly, the big
finitistic dimension $\mr{Findim}A$ is defined, allowing arbitrary
$A$-modules.

In 1960, Bass [\ref{Bs}] formulated two so-called finitistic
dimension conjectures. The first one asserts that
$\mr{findim}A=\mr{Findim}A$, while the second one claims that
$\mr{findim}A<\infty$.

It is known that the first finitistic dimension conjecture fails
in general [\ref{Zhi}] and the differences can even be arbitrarily
big [\ref{So}]. However, the second finitistic dimension
conjecture is still open. This conjecture is also related to many
other homological conjectures (e.g., the Gorenstein Symmetry
Conjecture, the Wakamatsu-tilting conjecture and the generalized
Nakayama conjecture) and attracts many algebraists, see for
instance [\ref{AhT}, \ref{Hp}, \ref{Xg}, \ref{Zho} etc.].

Only few classes of algebras was known to have finite finitistic
dimension, see for instance [\ref{GKK}, \ref{GZh}, \ref{Hp},
etc.]. In 1991, Auslander and Reiten [\ref{ARc}] proved that
$\mr{findim}A<\infty$ if $\mathcal{P}^{<\infty}$, the category of
all $A$-modules of finite projective dimension, is contravariantly
finite. However, Igusa-Smal$\phi$-Todorov [\ref{IST}] provided an
example of artin algebras over which $\mathcal{P}^{<\infty}$ is
not contravariantly finite.

In 2002, Angeleri-H\"{u}gel and Trlifaj [\ref{AhT}] obtained a
sufficient and necessary condition for the finiteness of
$\mr{findim}A$, using the theory of infinitely generated tilting
modules. They proved that $\mr{findim}A$ is finite if and only if
there exists a special tilting $A$-module (maybe infinitely
generated).

In the same year, Igusa and Todorov [\ref{ITf}] presented a good
way to prove the (second) finitistic dimension conjecture. In
particular, they proved that $\mr{findim}A$ is finite provided
that the representation dimension $\mr{repdim}A$, another
important homological invariant introduced by Auslander
[{\ref{Au}], is not more than 3. Recall that
$\mr{repdim}A=\mr{inf}\{\mr{gd}(\mr{End}_AV)\mid V$ is a
generator-cogenerator$\}$, where $\mr{gd}$ denotes the global
dimension and $\mr{End}_AV$ denotes the endomorphism algebra of
$_AV$. It should be noted that, though the representation
dimension of an artin algebra is always finite [\ref{Iy}], there
is no upper bound for the representation dimension of an artin
algebra [\ref{Rq}].

Using the Igusa-Todorov functor defined in [\ref{ITf}], Xi
[\ref{Xf}, \ref{Xg}, \ref{Xe}] developed some new ideas to prove
the finiteness of finitistic dimension of some artin algebras. For
example, it was shown in [\ref{Xg}] that, if $A$ is a subalgebra
of $R$ such that $\mr{rad}A$ is a left ideal in $R$ (in the case
$R$ is called a left idealized extension of $A$), then
$\mr{findim}A$ is finite provided that $\mr{repdim}R$ is not more
than 3. In addition, it was proved that the finitistic dimension
conjecture holds for all artin algebras if and only if the
finiteness of $\mr{findim}R$ implies the finiteness of
$\mr{findim}A$, for any pair $A,R$ such that $R$ is a left
idealized extension of $A$ [\ref{Xf}].

\vskip 15pt

In this note, we will continue the study of the finitistic
dimension conjecture following the ideas comes from [\ref{ITf},
\ref{Xf}, \ref{Xg}, \ref{Xe}]. Recall that an Auslander generator
over the artin algebra $A$ is a generator-cogenerator $_AV$ such
that $\mr{repdim}A=\mr{gd}(\mr{End}_AV)$. Note that it is not
known if $\mr{findim}A^{op}$ is finite provided that
$\mr{findim}A$ is finite in general, where $A^{op}$ denotes the
opposite algebra of $A$.

Our main results state as follows.

\bg{Th}\label{Mth}%
Let $A$ be an artin algebra with $\mr{repdim}A\le 3$. Assume that
$_AV$ is an Auslander generator. Then both
$\mr{findim}(\mr{End}_AM)$ and $\mr{findim}(\mr{End}_AM)^{op}$ are
finite, whenever $M\in\mr{add}_AV$.%
\ed{Th}%

%

Let $E$ be an artin algebra. Recall the following well known
conjectures (see for instance [\ref{Hp}, \ref{MR}], etc.).

{\bf Gorenstein Symmetry Conjecture}\ \  $\mr{id}(_EE)<\infty$ if
and only if $\mr{id}(E_E)<\infty$, where id denotes the injective
dimension.

{\bf Wakamatsu-tilting Conjecture}\ \  Let $_E\omega$ be a
Wakamatsu-tilting module. (1) If $\mr{pd}_E\omega<\infty$, then
$\omega$ is tilting. (2) If $\mr{id}_E\omega<\infty$, then
$\omega$ is cotilting.

{\bf Generalized Nakayama Conjecture}\ \  Each indecomposable
injective $E$-module occurs as a direct summand in the minimal
injective resolution of $_EE$.

Note that the Gorenstein symmetry conjecture and the Generalized
Nakayama conjecture are special cases of the second
Wakamatsu-tilting conjecture. Moreover, if the finitistic
dimension conjecture for $E$ and $E^{op}$, then all the above
conjectures hold.

\bg{Cor}\label{Corrd}%
Let $A$ be an artin algebra with $\mr{repdim}A\le 3$. Assume that
$_AV$ is an Auslander generator. Then the Gorenstein Symmetry
Conjecture, the Wakamatsu-tilting conjecture and the generalized
Nakayama conjecture hold for $\mr{End}_AM$.
\ed{Cor}%

%

As a special case of and Theorem \ref{Mth} and Corollary
\ref{Corrd}, we obtain the following result.

\bg{Cor}\label{Corfp}%
If $A=\mr{End}_{\Lambda}M$ for some module $M$ over a
representation-finite algebra $\Lambda$, then both $\mr{findim}A$
and $\mr{findim}A^{op}$ are finite.  In particular, the Gorenstein
Symmetry Conjecture, the Wakamatsu-tilting conjecture and the
generalized Nakayama conjecture hold for $A$.
\ed{Cor}%


We do not know whether or not every artin algebra can be realized
as an endomorphism algebra of some module over a
representation-finite algebra. If it is the case, then the above
result indicates that the finitistic dimension conjecture holds
for all artin algebras.

However,  it is well known that every artin algebra can be
realized as an endomorphism algebra of a projective and injective
module over a quasi-hereditary algebra. Thus the following result
shows that, in particular, if every quasi-hereditary algebra has a
left idealized extension which is a monomial algebra or an algebra
whose representation dimension is not more than 3, then  the
finitistic dimension conjecture holds.

\bg{Pro}\label{lid}%
Let $R$ be a left idealized extension of $A$. If $\mr{repdim}R\le
3$ or $\Omega^2_R(R\mr{-mod})$ is of finite type (for example,
$\mr{gd}R\le 2$), then $\mr{findim}(\mr{End}_AM)$ is finite, for
any projective $A$-module $M$.
\ed{Pro}

%
\vskip 30pt
\section{Finitistic dimension of endomorphism algebras}

\vskip 15pt


\hskip 15pt

Let $A$ be an artin algebra. We denoted by $A\mr{-mod}$ the
category of all $A$-modules. If $M$ is an $A$-module, we use
$\mr{pd}_AM$ to denote the projective dimension of $_AM$ and use
$\Omega^i_AM$ to denote the $i$-th syzygy of $M$. Throughout the
note, $\mathbf{D}$ denotes the usual duality functor between
$A\mr{-mod}$ and $A^{op}\mr{-mod}$.

The following lemma was well-known, see for instance [\ref{Xe}].

\bg{Lem}\label{rdlem}%
Let $A$ be an artin algebra and let $V$ be a generator-cogenerator
in $A\mr{-mod}$. The following are equivalent for a nonnegative
integer $n$.

(1) $\mr{gd}(\mr{End}_AV)\le n+2$.

(2) For any $X\in A\mr{-mod}$, there is an exact sequence $0\to
V_n\to \cdots\to V_1\to V_0\to X\to 0$ with each
$V_i\in\mr{add}_AV$, such that the corresponding sequence induced
by the functor $\mr{Hom}_A(V,-)$ is also exact.
\ed{Lem}%


The following lemma collects some important properties of the
Igusa-Todorov functor introduced in [\ref{ITf}].

\bg{Lem}\label{fclem}
For any artin algebra $A$, there is a functor $\Psi$ which is
defined on the objects of $A\mr{-mod}$ and takes nonnegative
integers as values, such that

$(1)$ $\Psi(M)=\mr{pd}_AM$ provided that $\mr{pd}_AM<\infty$.

$(2)$ $\Psi(X)\le \Psi(Y)$ whenever $\mr{add}_AX\subseteq
\mr{add}_AY$. The equation holds in case $\mr{add}_AX=
\mr{add}_AY$.

$(3)$ If $0\to X\to Y\to Z\to 0$ is an exact sequence in
$A\mr{-mod}$ with $\mr{pd}_AZ<\infty$, then $\mr{pd}_AZ\le
\Psi(X\oplus Y)+1$.
\ed{Lem}%


Let $A$ be an artin algebra and $M\in A\mr{-mod}$ with
$E=\mr{End}_AM$. Then $M$ is also a right $E$-module. It is well
known that $(M\otimes_E-, \mr{Hom}_A(M,-))$ is a pair of adjoint
functors and that, for any $E$-module $Y$, there is a canonical
homomorphism $\sigma_{Y}: Y\to \mr{Hom}_A(M,M\otimes_EY)$ definied
by $n\to [t\to t\otimes n]$. It is easy to see that $\sigma_Y$ is
an isomorphism provided that $Y$ is a projective $E$-module.

The following lemma is elementary but essential to prove our
results.

\bg{Lem}\label{Mlem}%
Let $M\in A\mr{-mod}$ and $E=\mr{End}_AM$. Then, for any $X\in
E\mr{-mod}$, $\Omega^2_EX\simeq \mr{Hom}_A(M,Y)$ for some $Y\in
A\mr{-mod}$.
\ed{Lem}%

\Pf. Obviously, we have an exact sequence

$$0\to \Omega^2_EX\to E_1\to E_0\to X\to 0$$

\noindent with $E_0, E_1\in E\mr{-mod}$ projective. Applying the
functor $M\otimes_E-$, we obtain an induced exact sequence

$$0\to Y\to M\otimes_E E_1\to M\otimes_E E_0\to
M\otimes_E X\to 0,$$

\noindent for some $Y\in A\mr{-mod}$. Now applying the functor
$\mr{Hom}_A(M,-)$, we further have an induced exact sequence

$$0\to \mr{Hom}_A(M,Y)\to \mr{Hom}_A(M,M\otimes_E
E_1)\to \mr{Hom}_A(M,M\otimes_E E_0).$$

Moreover, there is the following commutative diagram

\vskip 15pt

\setlength{\unitlength}{0.09in}
 \begin{picture}(60,8)

 \put(-1,1){\makebox(0,0)[c]{$0$}}
                             \put(0,1){\vector(1,0){2}}
 \put(9,1){\makebox(0,0)[c]{$\mr{Hom}_A(M,Y)$}}
                             \put(16,1){\vector(1,0){2}}
 \put(27,1){\makebox(0,0)[c]{$\mr{Hom}_A(M,M\otimes_E
E_1)$}}
                             \put(36,1){\vector(1,0){2}}
 \put(48,1){\makebox(0,0)[c]{$\mr{Hom}_A(M,M\otimes_E E_0).$}}

                 \put(10,4){\vector(0,-1){2}}
 \put(11,3.5){\makebox(0,0)[c]{$\phi$}}
                 \put(24,4){\vector(0,-1){2}}
 \put(26,3.5){\makebox(0,0)[c]{$\sigma_{E_1}$}}
                \put(48,4){\vector(0,-1){2}}
 \put(50,3.5){\makebox(0,0)[c]{$\sigma_{E_0}$}}

 \put(-1,6){\makebox(0,0)[c]{$0$}}
                             \put(0,6){\vector(1,0){2}}
 \put(10,6){\makebox(0,0)[c]{$\Omega^2_EX$}}
                             \put(16,6){\vector(1,0){2}}
 \put(24,6){\makebox(0,0)[c]{$E_1$}}
                             \put(36,6){\vector(1,0){2}}
 \put(48,6){\makebox(0,0)[c]{$E_0$}}

\end{picture}

\vskip 15pt

Since $E=\mr{End}_AM$ and $E_0, E_1\in \mr{add}_EE$, the canonical
homomorphisms $\sigma_{E_0}$ and $\sigma_{E_1}$ are isomorphisms.
It follows that $\Omega^2_EX\simeq \mr{Hom}_A(M,Y)$.

\hfill\#

\vskip 15pt
%

\bg{Th}\label{Mthm}%
Let $A$ be an artin algebra and $_AV$ be a generator-cogenerator
in $A\mr{-mod}$ such that $\mr{gd}(\mr{End}_AV)\le 3$. Then
$\mr{findim}E$ is finite, where $E=\mr{End}_AM$ for some
$M\in\mr{add}_AV$.
\ed{Th}

\Pf.  Suppose that $X\in E\mr{-mod}$ and that $\mr{pd}_EX<\infty$.
Then $\mr{pd}_E(\Omega^2_EX)<\infty$. Moreover, $\Omega^2_EX\simeq
\mr{Hom}_A(M,Y)$ for some $Y\in A\mr{-mod}$, by Lemma \ref{Mlem}.
Since $_AV$ is a generator and $\mr{gd}(\mr{End}_AV)\le 3$, by
Lemma \ref{rdlem} we obtain an exact sequence

$$0\to V_1\to V_0\to
Y\to 0\ (\dag)$$

\noindent with $V_0, V_1\in\mr{add}_AV$ such that the
corresponding sequence induced by the functor $\mr{Hom}_A(V,-)$ is
also exact. Note that $M\in\mr{add}_AV$, so the sequence ($\dag$)
also stays exact under the functor $\mr{Hom}_A(M,-)$. Thus, we
have the following exact sequence in $E\mr{-mod}$:

\centerline{$0\to \mr{Hom}_A(M,V_1)\to \mr{Hom}_A(M,V_0)\to
\mr{Hom}_A(M,Y)\to 0$.}

Now by Lemma \ref{fclem}, we deduce that

\hskip 115pt $\mr{pd}_EX\le \mr{pd}_E(\Omega^2_EX)+2$

\hskip 150pt $=\mr{pd}_E(\mr{Hom}_A(M,Y))+2$

\hskip 150pt
$\le\Psi(\mr{Hom}_A(M,V_0)\oplus\mr{Hom}_A(M,V_1))+1+2$

\hskip 150pt $\le\Psi(\mr{Hom}_A(M,V))+1+2<\infty$.

\hfill\#

\vskip 15pt

It is not known if $\mr{findim}A^{op}$ is finite provided that
$\mr{findim}A$ is finite in general. For an artin algebra $A$ with
$\mr{repdim}A\le 3$, it is known that $\mr{findim}A$ is finite
[\ref{Iy}] and since $\mr{repdim}A=\mr{repdim}A^{op}$,
$\mr{findim}A^{op}$ is also finite in the case. The following
corollary is also a generalization of this fact.

\bg{Cor}\label{Corrd2}%
Let $A$ be an artin algebra with $\mr{repdim}A\le 3$. Assume that
$_AV$ is an Auslander generator. If $M\in\mr{add}_AV$,  then both
$\mr{findim}(\mr{End}_AM)$ and $\mr{findim}(\mr{End}_AM)^{op}$ are
finite. In particular, the result holds for any
$_AM\in\mr{add}_A(A\oplus\mathbf{D}A)$.
\ed{Cor}%

\Pf. If $_AV$ is an Auslander generator in $A\mr{-mod}$, then
$_{A^{op}}\mathbf{D}V$ is an Auslander generator in
$A^{op}\mr{-mod}$. Moreover, if $M\in\mr{add}_AV$, then
$\mathbf{D}M\in\mr{add}_{A^{op}}\mathbf{D}V$. Obviously,
$(\mr{End}_AM)^{op}\simeq \mr{End}_{A^{op}}\mathbf{D}M$. Now the
conclusion follows from Theorem \ref{Mthm}.

\hfill\#
%
%
%
%
%
%

\vskip 15pt

Since all representation-finite algebras have representation
dimension 2, we obtain the following corollary.

\bg{Cor}\label{Corfp2}%
If $A=\mr{End}_{\Lambda}M$ for some representation-finite algebra
$\Lambda$ and some $M\in {\Lambda}\mr{-mod}$, then both
$\mr{findim}A$ and $\mr{findim}A^{op}$ are finite.
\ed{Cor}%

\Pf. Since $\Lambda$ is of finite representation type, the module
$V=\oplus_{i=1}^nV_i$ is an Auslander generator, where $\{V_i,
i=1,\cdots, n\}$ is the representation set of all indecomposable
$\Lambda$-modules. Now the conclusion follows from Corollary
\ref{Corrd2}.

\hfill\#

\vskip 15pt

The above result suggests a strong connection between the
finitistic dimension conjecture and the following endomorphism
algebra realization problem of artin algebras.

\vskip 15pt

\noindent {\bf Problem 1} Can all artin algebras be realized as
endomorphism algebras of modules over representation-finite
algebras? If not, what algebra has such a realization?

\vskip 15pt

Obviously, the affirmative answer to the first question implies
the finitistic dimension conjecture holds, by Corollary
\ref{Corfp2}.

In contrast, it is well known that every artin algebra can be
realized as an endomorphism algebra of a projective and injective
module over a quasi-hereditary algebra.

Let $A, R$ be both Artin algebras. Following [\ref{Xf}], we say
that $R$ is a left idealized extension of $A$ if $A\subseteq R$
has the same identity and $\mr{rad}A$ is a left ideal in $R$.


\bg{Pro}\label{lid2}%
Let $R$ be a left idealized extension of $A$. If $\mr{repdim}R\le
3$ or $\Omega^2_R(R\mr{-mod})$ is of finite type (for example,
$\mr{gd}R\le 2$), then $\mr{findim}E$ is finite, where
$E=\mr{End}_AM$ for some projective $A$-module $M$.
\ed{Pro}

\Pf.  Suppose that $X\in E\mr{-mod}$ and that $\mr{pd}_EX<\infty$.
Then $\mr{pd}_E(\Omega^2_EX)<\infty$. Moreover, $\Omega^2_EX\simeq
\mr{Hom}_A(M,Y)$ for some $Y\in A\mr{-mod}$, by Lemma \ref{Mlem}.
Since $_AM$ is projective, the proof of Lemma \ref{Mlem} in fact
indicates that $Y\in\Omega^2_A(A\mr{-mod})$. Then $Y\simeq
\Omega_RZ\oplus Q$ as $R$-modules for some $Z,Q\in R\mr{-mod}$
with $Q$ projective, by [\ref{Xf}, Erratum]. If $\mr{repdim}R\le
3$ or $\Omega^2_R(R\mr{-mod})$ is of finite type, we can obtain an
exact sequence of $R$-modules

\centerline{$0\to V_1\to V_0\to Y\to 0$ ($\ddag$)}

\noindent with $V_0, V_1\in\mr{add}_RV$. Here $_RV$ is an
Auslander generator in case $\mr{repdim}R\le 3$ or $V=R\oplus N$
in case $\Omega^2_R(R\mr{-mod})$ is of finite type, where $N$ is
the direct sum of nonisomorphic indecomposable $R$-modules
appeared in $\Omega^2_R(R\mr{-mod})$. The exact sequence ($\ddag$)
then restricts to an exact sequence in $A\mr{-mod}$. Since $_AM$
is projective, we have the following exact sequence in
$E\mr{-mod}$:

\centerline{$0\to \mr{Hom}_A(M,V_1)\to \mr{Hom}_A(M,V_0)\to
\mr{Hom}_A(M,Y)\to 0$.}

Now by Lemma \ref{fclem}, we deduce that

\hskip 115pt $\mr{pd}_EX\le \mr{pd}_E(\Omega^2_EX)+2$

\hskip 150pt $=\mr{pd}_E(\mr{Hom}_A(M,Y))+2$

\hskip 150pt
$\le\Psi(\mr{Hom}_A(M,V_0)\oplus\mr{Hom}_A(M,V_1))+1+2$

\hskip 150pt $\le\Psi(\mr{Hom}_A(M,V))+1+2<\infty$.

\hfill\#

\vskip 15pt

The above result has the following corollary which contains
[\ref{Xiab}, Corollary 2.15], where the result was proved under an
additional condition.

\bg{Cor}\label{mon}%
If $A$ has a left idealized extension which is a monomial algebra,
then $\mr{findim}E$ is finite, where $E=\mr{End}_AM$ for some
projective $A$-module $M$.
\ed{Cor}

\Pf. It was noted that $\Omega^2_R(R\mr{-mod})$ is of finite type
in case $R$ is a finite dimensional monomial relation algebra, see
[\ref{GZrp}, Example 2.3(a)] (also [\ref{Zhi}, Theorem A]). Thus
the conclusion follows from Proposition \ref{lid2}.

\hfill\#

\vskip 15pt

\noindent {\bf Problem 2} Does every quasi-hereditary algebra have
a left idealized extension which is a monomial algebra or an
algebra whose representation dimension is not more than 3?

\vskip 15pt

Of course, the affirmative answer to the problem also implies the
finitistic dimension conjecture holds.


{\small

}

\end{document}